\documentclass{amsart}
\usepackage{amssymb}

\DeclareFontEncoding{OT2}{}{} 
  \newcommand{\textcyr}[1]{%
    {\fontencoding{OT2}\fontfamily{wncyr}\fontseries{m}\fontshape{n}%
     \selectfont #1}}
\newcommand{\Sha}{{\mbox{\textcyr{Sh}}}}

\title{Curves Over Global Fields Violating The Hasse Principle} 
\author{Pete L. Clark}
\thanks{The author is partially supported by National Science Foundation grant DMS-0701771}
\address{Department of Mathematics, \\
                       University of Georgia, \\
                       Athens, GA 30602-7403, USA.}
\email{pete@math.uga.edu}

\begin{document}
\newtheorem{lemma}{Lemma}
\newtheorem{prop}[lemma]{Proposition}
\newtheorem{cor}[lemma]{Corollary}
\newtheorem{thm}[lemma]{Theorem}
\newtheorem{ques}{Question}
\newtheorem{quest}[lemma]{Question}
\newtheorem{conj}{Conjecture}
\newtheorem{fact}[lemma]{Fact}
\newtheorem*{mainthm}{Main Theorem}
\newtheorem{obs}[lemma]{Observation}
\newtheorem{hint}{Hint}
\maketitle

\newcommand{\pp}{\mathfrak{p}}
\newcommand{\DD}{\mathcal{D}}
\newcommand{\F}{\ensuremath{\mathbb F}}
\newcommand{\Fp}{\ensuremath{\F_p}}
\newcommand{\Fl}{\ensuremath{\F_l}}
\newcommand{\Fpbar}{\overline{\Fp}}
\newcommand{\Fq}{\ensuremath{\F_q}}
\newcommand{\PP}{\mathbb{P}}
\newcommand{\PPone}{\mathfrak{p}_1}
\newcommand{\PPtwo}{\mathfrak{p}_2}
\newcommand{\PPonebar}{\overline{\PPone}}
\newcommand{\N}{\ensuremath{\mathbb N}}
\newcommand{\Q}{\ensuremath{\mathbb Q}}
\newcommand{\Qbar}{\overline{\Q}}
\newcommand{\R}{\ensuremath{\mathbb R}}
\newcommand{\Z}{\ensuremath{\mathbb Z}}
\newcommand{\SSS}{\ensuremath{\mathcal{S}}}
\newcommand{\Rn}{\ensuremath{\mathbb R^n}}
\newcommand{\Ri}{\ensuremath{\R^\infty}}
\newcommand{\C}{\ensuremath{\mathbb C}}
\newcommand{\Cn}{\ensuremath{\mathbb C^n}}
\newcommand{\Ci}{\ensuremath{\C^\infty}}
\newcommand{\U}{\ensuremath{\mathcal U}}
\newcommand{\gn}{\ensuremath{\gamma^n}}
\newcommand{\ra}{\ensuremath{\rightarrow}}
\newcommand{\fhat}{\ensuremath{\hat{f}}}
\newcommand{\ghat}{\ensuremath{\hat{g}}}
\newcommand{\hhat}{\ensuremath{\hat{h}}}
\newcommand{\covui}{\ensuremath{\{U_i\}}}
\newcommand{\covvi}{\ensuremath{\{V_i\}}}
\newcommand{\covwi}{\ensuremath{\{W_i\}}}
\newcommand{\Gt}{\ensuremath{\tilde{G}}}
\newcommand{\gt}{\ensuremath{\tilde{\gamma}}}
\newcommand{\Gtn}{\ensuremath{\tilde{G_n}}}
\newcommand{\gtn}{\ensuremath{\tilde{\gamma_n}}}
\newcommand{\gnt}{\ensuremath{\gtn}}
\newcommand{\Gnt}{\ensuremath{\Gtn}}
\newcommand{\Cpi}{\ensuremath{\C P^\infty}}
\newcommand{\Cpn}{\ensuremath{\C P^n}}
\newcommand{\lla}{\ensuremath{\longleftarrow}}
\newcommand{\lra}{\ensuremath{\longrightarrow}}
\newcommand{\Rno}{\ensuremath{\Rn_0}}
\newcommand{\dlra}{\ensuremath{\stackrel{\delta}{\lra}}}
\newcommand{\pii}{\ensuremath{\pi^{-1}}}
\newcommand{\la}{\ensuremath{\leftarrow}}
\newcommand{\gonem}{\ensuremath{\gamma_1^m}}
\newcommand{\gtwon}{\ensuremath{\gamma_2^n}}
\newcommand{\omegabar}{\ensuremath{\overline{\omega}}}
\newcommand{\dlim}{\underset{\lra}{\lim}}
\newcommand{\ilim}{\operatorname{\underset{\lla}{\lim}}}
\newcommand{\Hom}{\operatorname{Hom}}
\newcommand{\Ext}{\operatorname{Ext}}
\newcommand{\Part}{\operatorname{Part}}
\newcommand{\Ker}{\operatorname{Ker}}
\newcommand{\im}{\operatorname{im}}
\newcommand{\ord}{\operatorname{ord}}
\newcommand{\unr}{\operatorname{unr}}
\newcommand{\B}{\ensuremath{\mathcal B}}
\newcommand{\Ocr}{\ensuremath{\Omega_*}}
\newcommand{\Rcr}{\ensuremath{\Ocr \otimes \Q}}
\newcommand{\Cptwok}{\ensuremath{\C P^{2k}}}
\newcommand{\CC}{\ensuremath{\mathcal C}}
\newcommand{\gtkp}{\ensuremath{\tilde{\gamma^k_p}}}
\newcommand{\gtkn}{\ensuremath{\tilde{\gamma^k_m}}}
\newcommand{\QQ}{\ensuremath{\mathcal Q}}
\newcommand{\I}{\ensuremath{\mathcal I}}
\newcommand{\sbar}{\ensuremath{\overline{s}}}
\newcommand{\Kn}{\ensuremath{\overline{K_n}^\times}}
\newcommand{\tame}{\operatorname{tame}}
\newcommand{\Qpt}{\ensuremath{\Q_p^{\tame}}}
\newcommand{\Qpu}{\ensuremath{\Q_p^{\unr}}}
\newcommand{\scrT}{\ensuremath{\mathfrak{T}}}
\newcommand{\That}{\ensuremath{\hat{\mathfrak{T}}}}
\newcommand{\Gal}{\operatorname{Gal}}
\newcommand{\Aut}{\operatorname{Aut}}
\newcommand{\tors}{\operatorname{tors}}
\newcommand{\Zhat}{\hat{\Z}}
\newcommand{\linf}{\ensuremath{l_\infty}}
\newcommand{\Lie}{\operatorname{Lie}}
\newcommand{\GL}{\operatorname{GL}}
\newcommand{\End}{\operatorname{End}}
\newcommand{\aone}{\ensuremath{(a_1,\ldots,a_k)}}
\newcommand{\raone}{\ensuremath{r(a_1,\ldots,a_k,N)}}
\newcommand{\rtwoplus}{\ensuremath{\R^{2  +}}}
\newcommand{\rkplus}{\ensuremath{\R^{k +}}}
\newcommand{\length}{\operatorname{length}}
\newcommand{\Vol}{\operatorname{Vol}}
\newcommand{\cross}{\operatorname{cross}}
\newcommand{\GoN}{\Gamma_0(N)}
\newcommand{\GeN}{\Gamma_1(N)}
\newcommand{\GAG}{\Gamma \alpha \Gamma}
\newcommand{\GBG}{\Gamma \beta \Gamma}
\newcommand{\HGD}{H(\Gamma,\Delta)}
\newcommand{\Ga}{\mathbb{G}_a}
\newcommand{\Div}{\operatorname{Div}}
\newcommand{\Divo}{\Div_0}
\newcommand{\Hstar}{\cal{H}^*}
\newcommand{\txon}{\tilde{X}_0(N)}
\newcommand{\sep}{\operatorname{sep}}
\newcommand{\notp}{\not{p}}
\newcommand{\Aonek}{\mathbb{A}^1/k}
\newcommand{\Wa}{W_a/\mathbb{F}_p}
\newcommand{\Spec}{\operatorname{Spec}}

\newcommand{\abcd}{\left[ \begin{array}{cc}
a & b \\
c & d
\end{array} \right]}

\newcommand{\abod}{\left[ \begin{array}{cc}
a & b \\
0 & d
\end{array} \right]}

\newcommand{\unipmatrix}{\left[ \begin{array}{cc}
1 & b \\
0 & 1
\end{array} \right]}

\newcommand{\matrixeoop}{\left[ \begin{array}{cc}
1 & 0 \\
0 & p
\end{array} \right]}

\newcommand{\w}{\omega}
\newcommand{\Qpi}{\ensuremath{\Q(\pi)}}
\newcommand{\Qpin}{\Q(\pi^n)}
\newcommand{\pibar}{\overline{\pi}}
\newcommand{\pbar}{\overline{p}}
\newcommand{\lcm}{\operatorname{lcm}}
\newcommand{\trace}{\operatorname{trace}}
\newcommand{\OKv}{\mathcal{O}_{K_v}}
\newcommand{\Abarv}{\tilde{A}_v}
\newcommand{\kbar}{\overline{k}}
\newcommand{\Kbar}{\overline{K}}
\newcommand{\pl}{\rho_l}
\newcommand{\plt}{\tilde{\pl}}
\newcommand{\plo}{\pl^0}
\newcommand{\Du}{\underline{D}}
\newcommand{\A}{\mathbb{A}}
\newcommand{\D}{\underline{D}}
\newcommand{\op}{\operatorname{op}}
\newcommand{\Glt}{\tilde{G_l}}
\newcommand{\gl}{\mathfrak{g}_l}
\newcommand{\gltwo}{\mathfrak{gl}_2}
\newcommand{\sltwo}{\mathfrak{sl}_2}
\newcommand{\h}{\mathfrak{h}}
\newcommand{\tA}{\tilde{A}}
\newcommand{\sss}{\operatorname{ss}}
\newcommand{\X}{\Chi}
\newcommand{\ecyc}{\epsilon_{\operatorname{cyc}}}
\newcommand{\hatAl}{\hat{A}[l]}
\newcommand{\sA}{\mathcal{A}}
\newcommand{\sAt}{\overline{\sA}}
\newcommand{\OO}{\mathcal{O}}
\newcommand{\OOB}{\OO_B}
\newcommand{\Flbar}{\overline{\F_l}}
\newcommand{\Vbt}{\widetilde{V_B}}
\newcommand{\XX}{\mathcal{X}}
\newcommand{\GbN}{\Gamma_\bullet(N)}
\newcommand{\Gm}{\mathbb{G}_m}
\newcommand{\Pic}{\operatorname{Pic}}
\newcommand{\FPic}{\textbf{Pic}}
\newcommand{\solv}{\operatorname{solv}}
\newcommand{\Hplus}{\mathcal{H}^+}
\newcommand{\Hminus}{\mathcal{H}^-}
\newcommand{\HH}{\mathcal{H}}
\newcommand{\Alb}{\operatorname{Alb}}
\newcommand{\FAlb}{\mathbf{Alb}}
\newcommand{\gk}{\mathfrak{g}_k}
\newcommand{\car}{\operatorname{char}}
\newcommand{\Br}{\operatorname{Br}}
\newcommand{\gK}{\mathfrak{g}_K}
\newcommand{\coker}{\operatorname{coker}}
\newcommand{\red}{\operatorname{red}}
\newcommand{\CAY}{\operatorname{Cay}}
\newcommand{\ns}{\operatorname{ns}}
\newcommand{\xx}{\mathbf{x}}
\newcommand{\yy}{\mathbf{y}}
\newcommand{\E}{\mathbb{E}}
\newcommand{\rad}{\operatorname{rad}}
\newcommand{\Top}{\operatorname{Top}}
\newcommand{\Map}{\operatorname{Map}}
\newcommand{\Li}{\operatorname{Li}}
\newcommand{\Res}{\operatorname{Res}}
\newcommand{\Sel}{\operatorname{Sel}}
\renewcommand{\sep}{\operatorname{sep}}
\renewcommand{\Gal}{\operatorname{Gal}}
\renewcommand{\Hom}{\operatorname{Hom}}
\renewcommand{\GL}{\operatorname{GL}}
\newcommand{\SL}{\operatorname{SL}}
\renewcommand{\div}{\operatorname{div}}
\newcommand{\nn}{\mathfrak{n}}




\begin{abstract}
In response to a question of B. Poonen, we exhibit for each global field $k$ an algebraic curve over $k$ 
which violates the Hasse Principle.  In fact we can find such examples among Atkin-Lehner twists of certain 
elliptic modular curves and -- in positive characteristic -- Drinfeld modular curves.  Our main tool is a 
refinement of the ``Twist Anti-Hasse Principle'' (TAHP).  We also use TAHP to construct further Hasse 
Principle violations, for instance among curves over any number field of any given genus $g \geq 2$.
\end{abstract}  

\section*{Notation}
\noindent
We let $k_0$ denote either $\Q$ or the rational function field $\F_p(t)$ for a prime number $p$.  By a 
\textbf{global field} we mean a finite separable extension $k$ of $k_0$. 
\\ \\
If $v$ is a place of $k$, we denote the completion by $k_v$.  If $v$ is non-Archimedean, we denote by 
$R_v$ the ring of integers of $k_v$, $\mathfrak{m}_v$ the maximal ideal of $R_v$ and $\F_v = R_v/\mathfrak{m}_v$ 
the residue field.  Let $\tilde{\Sigma}_k$ denote the set of all places of $k$ and $\Sigma_k$ the set of all 
non-Archimedean places of $k$.  We write $\mathbb{A}_k$ for the ad\`ele ring over $k$.
\\ \\
By a \textbf{nice} variety $V_{/k}$, we mean an algebraic $k$-variety which is smooth, projective and 
geometrically integral.  A nice curve $C_{/k}$ is a nice algebraic variety of dimension one.
\\ \\
Let $k$ be a field.  The \textbf{characteristic exponent} of $k$ is $1$ if $\car(k) = 0$ and otherwise is 
equal to $\car(k)$.  
\\ \\
Let $p$ be either $1$ or a prime number.  Let $k$ be a global field of characteristic 
exponent $p$, and let $M,N$ be odd prime 
numbers such that $M \geq 4$ and $M,N,p$ are pairwise coprime.  Then the moduli problem which associates to a 
$\Z[\frac{1}{MN}]$-scheme $S$ the set of isomorphism classes of triples $(E,C,P)$, where $E_{/S}$ is an elliptic curve, $C \subset_S E[N]$ is a locally 
cyclic subgroup scheme of order $N$ and $P$ is an $S$-rational point of order $M$ is representable 
by a smooth, projective relative curve $X(N,M)_{/\Z[\frac{1}{MN}]}$.  In particular, the modular curve 
$X(N,M)$ has a model over $\Q$ and over $\F_p$, and hence also over $k$.

\section{Introduction}
\noindent
In February of 2009, Bjorn Poonen wrote to the author with the following question:
\\ \\
\indent ``[D]o you know if this is a 
known result?  There is an algorithm that takes as input a number field k, and produces a curve over k 
that violates the Hasse principle.  (I just observed yesterday that this can be proved.)'' 
\\ \\
The author replied that he thought the result was interesting, was to his knowledge unknown, and that he might be able to 
prove it using some ``old tricks.''  But he did not pursue the matter further until seeing Poonen's 
preprint \cite{Poonen09}.  In March of 2009, the author wrote to Professor Poonen and sketched his proof.  
This sketch now appears in \cite{Poonen09}, along with a sketch of a third proof by 
J.-L. Colliot-Th\'el\`ene.    
\\ \\
It is the main goal of the present paper to present a detailed proof of the author's construction, extended to 
include the function field case:
\begin{mainthm}
\label{MAINTHM}
Let $k/k_0$ be a global field of characteristic exponent $p$ and degree $d = [k:k_0]$.   Then there exists an 
effectively computable constant $B = B(p,d)$ such that: for any odd prime numbers $M,\  N > B$ with $N \equiv -1 \pmod M$, there exists an infinite, effectively computable set $\{l_i\}$ of separable quadratic extensions 
of $k_0$ such that for all $i$, the twist $\mathcal{T}_i$ of $X(N,M)$ by the Atkin-Lehner involution $w_N$ and 
$l_i/k_0$ has points everywhere locally over $k_0$ and no $k$-rational points.  It follows that for all $i$, 
the base extension $\mathcal{T}_i \otimes_{k_0} k$ of $\mathcal{T}_i$ to $k$ violates the Hasse Principle over $k$.
\end{mainthm}
\noindent
We also effectively construct, for every global field $k$ of positive characteristic, an Atkin-Lehner twisted 
Drinfeld modular curve which violates the Hasse Principle over $k$: Theorem \ref{DRINFELDTHM}.
\\ \\
Thus our Main Theorem not only answers Poonen's question, but shows (we believe) that such violations are not 
a contrivance but arise already in arithmetically natural families of curves: the curves $\mathcal{T}_i$ 
parameterize certain elliptic $\Q$-curves endowed with a torsion point of order $M$.
\\ \\
As alluded to above, the ``old trick'' that we use to prove the Main Theorem is the \textbf{Twist Anti-Hasse Principle} (TAHP) of \cite{ClarkPRIME}.  The version of TAHP presented in 
\emph{loc. cit.} applies to curves over $\Q$.  In $\S 2$ we prove a version of TAHP which is different in 
three respects.  First, we work in the context of global fields of arbitrary characteristic.  Second, we
give additional conditions which suffice for the HP violation to persist over certain extensions of the ground 
field.  Third, we no longer emphasize ``prime twists'' in the statements, although such considerations continue 
to play a role in the proofs in the form of using quadratic extensions which are ramified at as few finite 
places as possible.  The TAHP is proved in $\S 2$, and in $\S 3$ it is applied, together with results of 
Merel and Poonen, to prove the Main Theorem.
\\ \\
One of the author's continuing research interests is to produce a large and varied supply of curves over global fields violating 
HP.  It is natural to ask for HP violations with prescribed arithmetic-geometric invariants. 

\begin{ques}
Can one find -- effectively, if possible -- \\
a) For each global field $k$, a genus one curve $C_{/k}$ violating HP? \\
b) For each global field $k$ and $g \geq 2$, a genus $g$ curve $C_{/k}$ violating HP? \\
c) As in c), but (for $g \geq 3$) nonhyperelliptic?  \\
d) As in c), but (for fixed $d$ and $g \gg d$) of gonality $d$? \\
e) As in c) but (for $g \geq 3$) but with no nontrivial automorphisms?  
\end{ques}
\noindent 
The TAHP can only produce HP violations in genus $g \geq 2$.  In $\S 4$, we use it 
to answer some of the above questions as they pertain to the $g \geq 2$ case.  In $\S 5$ we discuss some 
conjectures on HP violations in genus one and interrelationships between them.  In $\S 6$ we explain how it 
is much easier to construct higher-dimensional varieties $V$ over any global field $k$ violating HP and discuss 
the prospect of passing from a variety violating HP to a subvariety which still violates HP.  
\\ \\
Acknowledgments: Thanks to B. Poonen for suggesting the problem and for several helpful pointers to the literature, to E. Izadi for the proof of 
Lemma \ref{IZADILEMMA}, and to D. Swinarski and D. Lorenzini for helpful conversations.

\section{The Twist Anti-Hasse Principle Revisited}
\noindent
\subsection{The Twist anti-Hasse principle}

%
\begin{thm}(Twist Anti-Hasse Principle)
\label{TAHP}
Let $C_{/k}$ a nice curve and $\iota: C \ra C$ a $k$-rational involutory automorphism 
$(\iota^2 = 1)$.  Let $C/\iota$ be the quotient of $C$ by the group $\langle \iota \rangle$, 
so that there is a natural map $\Psi: C \ra C/\iota$ of degree $2$ and inducing a separable quadratic extension of 
function fields.  Suppose: \\
(i) $\{P \in C(k) \ | \ \iota(P) = P \} = \varnothing$. \\
(ii) $\{P \in C(\overline{k}) \ | \ \iota(P) = P\} \neq \varnothing$. \\
(iii) $C(\mathbb{A}_k) \neq \varnothing$. \\
(iv) $\# (C/\iota)(k) < \infty$.  \\
Then there exist infinitely many separable quadratic extensions $l/k$ such that the twisted curve 
$\mathcal{T}(C,\iota,l/k)$ violates the Hasse Principle over $k$.  We may take each $l/k$ to 
be the base change of a separable quadratic extension of the subfield $k_0$. 
\end{thm}
\begin{proof}
Let $l/k$ be a separable quadratic extension.  Then by $\mathcal{T}_l(C) := \mathcal{T}(C,\iota,l/k)$ we mean the $k$-variety 
obtained from the $l$-variety $C_{/l}$ by twisting the given $k$-structure by the cohomology class 
corresponding to $l$ in \[H^1(k,\langle \iota \rangle) \cong H^1(k,\Z/2 \Z) = \Hom(\Gal(k^{\sep}/k),\Z/2\Z). \] More concretely, 
letting $\sigma_l$ denote the nonidentity element of $\Gal(l/k)$, the twisted $\Gal(l/k)$-action on 
$C(l)$ is given by $P \mapsto \iota(\sigma_l P)$.  For each such $l$, we have natural set maps
\[\alpha_l: \mathcal{T}_l(C)(k) \hookrightarrow C(l), \]
\[\beta_l: C(l) \ra (C/\iota)(l). \]
 Put
\[S_l = (\beta_l \circ \alpha_l)(C_l(k)). \] 
Then $S_l\subset
(C/\iota)(k)$. Moreover, $(C/\iota)(k) = \bigcup_l S_l \cup
\Psi(C(k))$, and for $l \neq l'$, $P \in S_l \cap S_{l'}$
implies that $P \in C(l) \cap C(l') = C(k)$.  But $S_l \cap C(k)$ consists of $k$-rational
$\iota$-fixed points, which we have assumed in (i) do not exist,
so that for $l \neq l'$, $S_l \cap S_{l'} = \emptyset$. By (iv),
$(C/\iota)(k)$ is finite, and we conclude that the set of $l$ for
which $S_l \neq \emptyset$ is finite.  
\\ \\
It therefore suffices to find infinitely many separable quadratic extensions $l$ 
such that $\mathcal{T}_l(C)$ has points everywhere locally.
\\ \\
We now consider cases, depending upon whether the characteristic of $k$ is zero, positive and odd, or equal to $2$.  
\\ \\
Case 1: $k$ has characteristic $0$ (i.e., $k$ is a number field).  Then Kummer theory applies to all 
quadratic extensions: every quadratic extension $l/k$ is separable and of the form $l = k(\sqrt{d})$ 
for some $d \in k \setminus k^{\times 2}$.  Moreover $k(\sqrt{d}) = k(\sqrt{d'})$ iff 
$d \equiv d' \pmod {k^{\times 2}}$.  Let us in fact take $l = k(\sqrt{d})$ with $d \in \Q$ of $k$.  It will 
be notationally convenient to think of the choice of $d$ -- and hence $l = k(\sqrt{d})$ -- as fixed: we 
may then abbreviate 
\[\mathcal{T} C := \mathcal{T}_{k(\sqrt{d})}C. \]
To be sure, our main task is to demonstrate that there are infinitely many classes $d \in \Q^{\times}/\Q^{\times 2}$ 
such that $\mathcal{T} C$ has points everywhere locally.     
\\ \\
First we require $d > 0$.  This ensures that for every real place $v$ of $k$ (if any) we have 
$\mathcal{T}_l(C)(k_v) = C(k_v) \neq \varnothing$.  Henceforth we need only worry about the finite places.  
\\ \\
Let $R_d$ be the set of finite places that ramify in $l = k(\sqrt{d})$.  Let $M_1$ be a positive integer such that for all places $v$ with $\# \F_v > M_1$, 
$C$ extends to a smooth relative curve $C_{/R_v}$ and such that every smooth curve of genus $g(C)$ over a 
finite field $\F$ with at least $M_1$ elements has an $\F$-rational point.  We call a finite place $v$ 
\textbf{large} if $\# \F_v > M_1$ and $v \not \in R_d$; otherwise $v$ is \textbf{small}.  Note that for any 
fixed $d$, all but finitely many places are large.  
\\ \\
We claim that for any large place $v$, $\mathcal{T} C (k_v) \neq 
\varnothing$.  To see this, consider the minimal regular $R_v$-model for $\mathcal{T} C$.  Since $v$ is 
large and $\mathcal{T} C$ becomes isomorphic to $C$ after the base change $k \mapsto k(\sqrt{d})$, 
the minimal regular $R_v$-model for $\mathcal{T} C$ becomes smooth after an unramified base change.  However, 
smoothness is a geometric property and formation of the regular model commutes with unramified base change, 
so it must be the case that $\mathcal{T} C$ itself extends smoothly to $R_v$.  Moreover, by assumption on $v$, 
there is an $\F_v$-rational point on the special fiber, so by Hensel's Lemma 
$\mathcal{T} C(k_v) \neq \varnothing$.  
\\ \\
Now let $P$ be an $\iota$-fixed point, and put $K = k(P)$.  Suppose that we can choose $d$ such that all small 
places $v$ split completely in $K$.  Then since $\mathcal{T} C$ has $K$-rational points, it has $k_v$-rational 
points for all small places $v$, and hence it has points everywhere locally.
\\ \\
So it is enough to show the following: there exist infinitely many $d \in \Q^{\times}/Q^{\times 2}$ 
such that: $d > 0$, $k(d)/k$ is a quadratic extension, and for each finite place $v$
of $k$ such that either (i) $\# \F_v \leq M_1$ or (ii) $v$ ramifies in $k(d)/k$, we have that $v$ splits 
completely in $k(P)/k$.  Let $\zeta_4$ be a primitive $4$th root of unity.  By Cebotarev density, the set 
$\mathcal{P}$ of primes of $\Q$ which split completely in the finite separable extension $k(P,\zeta_4)/\Q$ has 
positive density.  More concretely, under the usual identification of finite places of $k_0$ with 
prime numbers, any element of $\mathcal{P}$ is a prime number $p \equiv 1 \pmod 4$, so that $k(\sqrt{p})/k$ 
ramifies only at primes of $k$ lying over $p$.  Thus taking $d$ to be any element of $\mathcal{P}$, we get that 
$\mathcal{T} C$ has points everywhere locally. \\ \\
Case 2: $k$ is a function field of odd characteristic $p$.  Much of Case 1 still applies: in particular, 
quadratic extensions are governed by Kummer theory.  The main difference is that the prime subfield is now 
$\F_p(t)$, so that its places correspond to irreducible polynomials $p_i(t) \in \F_p[t]$, together with the place 
corresponding to the point at infinity on the projective line.  Arguing as above, it suffices to find infinitely 
many squarefree polynomials $d(t) \in \F_p[t]$ such that the corresponding quadratic extension 
$l = \F_p(t)(\sqrt{d(t)})/\F_p(t))$ ramifies only at places of $k_0$ which split completely in $k(P)$.  But the 
ramified places of $l/k$ correspond to the irreducible factors of $d(t)$, together with the place at infinity if and only if $d(t)$ has 
odd degree.  By Cebotarev Density, there are infinitely many irreducible polynomials 
$p_i(t) \in \F_p[t]$ which split completely in $k(P)/k_0$.  So taking $d = p_i(t)$ works provided $p_i(t)$ has 
even degree.  Thus, if infinitely many $p_i$'s have even degree, we're done.  Otherwise, by passing to a subseqence 
we may assume that all the $p_i$'s have even degree and take $d = f_i f_j$ for $i \neq j$. 
\\ \\
Case 3: Suppose that $k$ has characteristic $2$.  Again most aspects of the proof go through, but in place of Kummer theory we must use Artin-Schreier 
theory, which turns out to be simpler.  By Cebotarev, there exist infinitely many irreducible polynomials 
$p_i(t) \in \F_2[t]$ which split completely in the separable extension $k(P)/\F_2(t)$.  The quadratic extension 
$l/k$ defined by the polynomial 
$X^2 + X = \frac{1}{p_i(t)}$
is then unramified at every place $v$ such that $v_p(\frac{1}{p_i(t)}) \geq 0$, i.e., at every place except 
the place corresponding to $p_i$.  This completes the proof. 
\end{proof}

\subsection {A corollary on base extension}

\begin{cor}
\label{TWISTAHPCOR}
Let $k/k_0$ be a global field.  Suppose that we have a curve $C_{/k}$ and a $k$-rational involution 
$\iota$ satisfying hypotheses (i) through 
(iv) of Theorem \ref{TAHP}.  Suppose further that $K/k$ is a finite separable 
field extension such that the following hypotheses all hold: \\
$(i)_K$: There are no $K$-rational $\iota$-fixed points; \\
$(iv)_K$: $\# (C/\iota)(K) < \infty$. \\ Then there are infinitely many separable quadratic extensions $l/k_0$ 
such that the twisted curve $\mathcal{T}(C,\iota,lk/k) \otimes_{k} K$ violates the 
Hasse Principle over $k$.
\end{cor}
\begin{proof}
Observe that hypotheses (ii) and (iii) are stable under base change: since they hold over $k$, \emph{a fortiori} 
they also hold over $K$.  We may therefore apply Theorem \ref{TAHP} to the pair $(C,\iota)_{/K}$.  For all 
but finitely many $l$'s, $Kl/K$ is a separable quadratic extension, and for such $l$ we have
\[\mathcal{T}(C_{/K},\iota_{/K},Kl/K) = \mathcal{T}(C,\iota,lk/k) \otimes_{k} K. \] 
\end{proof}
\subsection{Remarks on the satisfaction of the hypotheses} \textbf{} \\ \\ \noindent
Recall \cite[Remark 1.2]{ClarkPRIME}: if a pair $(C,\iota)_{/k}$ satisfies the hypotheses of Theorem 
\ref{TAHP} then the genus of $C$ is at least $2$ and the genus $g(C/\iota)$ of $C/\iota$ is at least one.  
\\ \\
Let us separately consider some cases:
\\ \\
$\bullet$ $C/\iota$ has genus one and $0 < \# (C/\iota)(k) < \infty$.  Then $C/k = E$ is an elliptic 
curve with finitely many $k$-rational points.  Under these circumstances we can construct many pairs
$(C,\iota)_{/k}$ satisfying the hypotheses of Theorem \ref{TAHP} and such that $C/\iota \cong E$: see
Theorem \ref{BIELLIPTICTHM} below.  However, I do not know how to systematically produce extensions $K/k$ 
such that the hypotheses of Corollary \ref{TAHP} are satisfied: this is the notorious problem of 
growth of the Mordell-Weil 
group under field extension. 
\\ \\
$\bullet$ $C/\iota$ has genus one and is a locally trivial but globally nontrivial torsor under its Jacobian 
elliptic curve $E$.  In this case $C/\iota$ itself already violates the HP over $k$.  Moreover, upon base-changing 
to any finite extension $K/k$ which is not divisible by the index of $C/\iota$, we retain a HP violation.
\\ \\
$\bullet$ $k$ is a number field and $C/\iota$ has genus at least $2$.  In this case, Faltings' theorem 
asserts that $\# C(K) < \infty$ for every finite extension $K/k$, i.e., hypothesis (iv)$_K$ of Corollary 
\ref{TWISTAHPCOR} holds for all extensions $K/k$.  Therefore in order to apply Corollary \ref{TWISTAHPCOR} 
to $(C,\iota)_{/k}$ and $K/k$, we need only verify $(i)_K$: that $\iota$ has no $K$-rational fixed points.  
But this is a very mild hypothesis: if $C$ has no $k$-rational $\iota$-fixed points, then the set of 
$\overline{k}$-fixed points breaks up into a finite set of $\Gal(\overline{k}/k)$-orbits, each of size 
greater than one.  This means that there 
is a finite set of nontrivial Galois extensions $m_1,\ldots,m_d/k$ such that for any finite extension 
$K/k$, (i)$_K$ holds iff for all $i$, $K \not \supset m_i$.  
\\ \indent
Suppose now that instead of a single pair $(C,\iota)_{/k}$ we have an infinite sequence $(C_n,\iota_n)_{/k}$ of such pairs, then we get an infinite 
sequence of finite sets of nontrivial Galois extensions: $(\{m_{n,1},\ldots,m_{n,d_n)} \})_{n=1}^{\infty}$.  Make 
the following additional hypothesis (H): for every finite extension $K/k$, there exists $n \in \Z^+$ 
such that $K \not \supset m_{n,i}$ for all $1 \leq i \leq d_n$.  Then (i)$_K$, (ii), (iii), (iv) 
hold for $(C_n,\iota_n)$ and $K/k$, so that we get infinitely many HP violations $\mathcal{T}(C_n,\iota_n,Kl/K)$ 
over $K$.  
\\ \\
$\bullet$ $k$ is a function field, $C/\iota$ has genus at least two and has transcendental moduli 
(``nonisotrivial'').  Then a theorem of P. Samuel \cite[Thm. 4]{Samuel} asserts that, as above, for any finite field extension 
$K/k$, $\# (C/\iota)(K) < \infty$, and the discussion is the same as in the previous case.
\\ \\
$\bullet$ $k$ is a function field, $C/\iota$ has genus at least two and has algebraic moduli (``isotrivial''): 
in other words, $(C/\iota)_{/\overline{k}}$ has a model over some finite field $\F_q$.  In this case the set $(C/\iota)(K)$ need not be 
finite, and indeed, after enlarging the base so that $C/\iota$ is defined over a finite subfield $\F_q$ of $k$, 
as soon as there exists a single point $P \in (C/\iota)(k) \setminus (C/\iota)(\F_q)$, there are infinitely many 
$k$-rational points: iterates of the $q$-power Frobenius map fix the curve but not the point $P$.  In this case we can 
show finiteness by arranging for a much stronger property to hold: writing $\F$ for the algebraic closure of 
$\F_p$ in $k$, we certainly have that $C(\F)$ is finite, so $C(k)$ is finite if $C(k) = C(\F)$. 

\section{HP-violations from Atkin-Lehner twists of modular curves} 

\subsection{Number field case} \textbf{} \\ \\
\noindent
We begin with the following setup, taken from \cite{ClarkPRIME}.  Let $N\in \Z^+$ be squarefree, 
let $X_0(N)_{/\Q}$ denote the modular curve with $\Gamma_0(N)$-level structure, and let $w_N$ denote the Atkin-Lehner 
involution, a $\Q$-rational involutory automorphism of $X_0(N)$.  Let $k/\Q$ be a number field.  We claim 
there exists an effectively computable $N_0$ such that for all squarefree $N > N_0$, 
$(X_0(N),w_N)$ satisfies hypotheses (i)$_k$, (ii), (iii), (iv)$_k$ of Corollary \ref{TWISTAHPCOR}, and therefore 
there are infinitely many quadratic twists of $X_0(N)$ which violate the Hasse Principle over $\Q$ and also 
over $k$.
\\ \indent
To see this, recall the following information from \cite{ClarkPRIME} and the references cited therein: there 
are always geometric $w_N$-fixed points -- i.e., hypothesis (ii) holds -- and the least degree $[\Q(P):\Q]$ of a $w_N$-fixed point is equal 
to the class number of $\Q(\sqrt{-N})$.  By a famous result of Heilbronn, the class number of $\Q(\sqrt{-N})$
tends to infinity with $N$, so that hypothesis (i)$_k$ holds for each number field $k$ and all $N \geq N_0(k)$, 
as above.  Further, the cusp at $\infty$ is a $\Q$-rational point on $X_0(N)$, so (iii) holds.  Moreover, 
the genus of $X_0^+(N) := X_0(N)/w_N$ is at least $2$ for all $N \geq 133$.  Therefore for $N \geq 
\max(N_0(k),133)$, there are infinitely many quadratic extensions $l/k$ such that $\mathcal{T}(X_0(N),w_N,kl/k)$ 
violates the Hasse Principle.  
\\ \\
In this construction, it is easy to make $N_0(k)$ explicit, but in order to get an explicit choice of $l$ we 
would need to explicitly determine all points of $X_0(N)$ which are defined over a quadratic extension of 
the number field $k$.  This lies beyond the current state of knowledge of rational points on $X_0(N)$.  However, 
we have a much better understanding of rational points on the coverings $X_1(N)$:
\begin{thm}(Merel \cite{Merel})
\label{MERELTHM}
For each positive integer $d$, there exists an effectively computable integer $M(d)$ such that if $P \in X_1(M)$ is any 
noncuspidal point with $[\Q(P):\Q] \leq d$, then $M \leq M(d)$.
\end{thm}
\noindent
Now let $N$ and $M$ be coprime squarefree positive integers with $M \geq 4$, and let $X(N,M)_{/\Q}$ be the modular curve with 
$\Gamma_0(N) \cap \Gamma_1(M)$-level structure.  As recalled in $\S 1$, $X(N,M)_{/\Q}$ is a fine moduli space 
for the moduli problem $(E,C,P)_{/S}$: here $S$ is a $\Q$-scheme, $E_{/S}$ is an elliptic curve, $C \subset_S E[N]$ is a cyclic order $N$-subgroup scheme, and $P \in E(S)$ is a 
point of order $M$.  The involution $w_N$ is defined on moduli as follows: 
\[w_N: (E,C,P) \mapsto (E/C,E[N]/C,\iota(P)), \]
where $E/C$ is the quotient elliptic curve and $\iota$ is the isogeny $E \ra E/C$.  The 
canonical map $X(N,M) \ra X_0(N)$ can be viewed on moduli as the forgetful map $(E,C,P) \mapsto (E,C)$, 
and the action of $w_N$ on $X_0(N)$ is the evident compatible one: $(E,C) \mapsto (E/C,E[N]/C)$.  It follows 
that the $w_N$-fixed points of $X(N,M)$ -- if any -- lie over the $w_N$-fixed points of $X_0(N)$.  Therefore, 
so long as $N$ is chosen sufficiently large with respect to the number field $k$ as above, since (i)$_k$ 
holds for $w_N$ on $X(N)$, \emph{a fortiori} it holds for $w_N$ on $X(N,M)$.  As for any modular curve 
corresponding to a congruence subgroup, the cusp at $\infty$ gives a rational point on $X(N,M)$, so (iii) holds.  
Similarly, since there is a natural map $X(N,M)/w_N \ra X_0(N)/w_N$, the curve $X(N,M)/w_N$ has genus at least 
$2$ for all $N \geq 133$ and therefore (iv)$_{/k}$ holds.  Indeed, because of Merel's theorem we can be more explicit.  
\\ \\
But first we need to verify (ii): that $X(N,M)$ has any $w_N$-fixed points at all.  For any squarefree $N > 2$, 
one set of $w_N$-fixed points on $X_0(N)$ corresponds to the set of distinct isomorphism classes of elliptic 
curves $E$ with $\Z[\sqrt{-N}]$-CM, and $C$ is the kernel of the endomorphism $\iota = [\sqrt{-N}]$.  (If 
$N \equiv 1 \pmod 4$, this is the entire set of $w_N$-fixed points.  Otherwise, there is another set corresponding 
to elliptic curves with CM by the maximal order of $\Q(\sqrt{-N})$.)  In order to get fixed points on $X(N,M)$, 
we wish to find an $M$-torsion point on a $\Z[\sqrt{-N}]$-CM elliptic curve with $\iota(P) = P$.  Observe that 
$\iota(P) = P$ implies $P = \iota^2(P) = -NP$, so this is only possible
if $(N+1)P = 0$, i.e., if $N \equiv -1 \pmod M$.  Conversely, if $N \equiv -1 \pmod M$, then the characteristic 
polynomial of $\iota$ acting on $E[M]$ is $t^2 + N = t^2 - 1 = (t+1)(t-1)$, so there is an eigenvalue of $1$, 
i.e., an $\iota$-fixed point. 
\\ \\
Now we put all the ingredients together: let $k$ be any number field; put $d = [k:\Q]$.  Let $M$ be an odd prime
number which is larger than Merel's bound $M(2d)$.  By Dirichlet's theorem 
on primes in arithmetic progressions, there exist infinitely many prime numbers $N$ such that $N \equiv -1 \pmod M$; choose one with $h(\Q(\sqrt{-N})) > d$.  Then the curve $X(N,M)_{/\Q}$ satisfies 
(i)$_k$, (ii), (iii), (iv) of Corollary \ref{TWISTAHPCOR}.  Then $X(N,M)$ has no nonscuspidal rational points in 
any quadratic extension of $k$, so that there is an infinite, effectively computable set of quadratic extensions 
$l/\Q$ such that $\mathcal{T}(X(N,M),w_N,kl/k)$ violates the Hasse Principle over $k$.
\\ \\
We can say a bit more about the quadratic fields $l/\Q$ which work in the above construction. $X(N,M)$ is a quotient of 
the modular curve $X(NM)_{/\Q(\zeta_{NM}}$ with full $NM$-level structure.  The cusps form a single orbit 
under the automorphism group of $X(NM)$.  Since the cusp at $\infty$ is $\Q$-rational, it follows that 
all of the cusps are rational over the cyclotomic field  $\Q(\zeta_{NM})$.  Since none of the cusps are 
$w_N$-fixed points, it follows that if $l$ is not contained in $k(\zeta_{NM})$ the cuspidal points will not 
be rational over $\mathcal{T}(X(N,M),w_N,kl/k)_{k(\zeta_{NM})}$: thus only an explicit finite set of quadratic 
fields $l/\Q$ must be excluded.  

\subsection{Function field case}
\textbf{} \\ \\
Suppose now that $p$ is a prime number and $k_0 = \F_p(t)$.  The situation here differs from the number field case in that we cannot use 
Theorem \ref{TAHP} to get Atkin-Lehner twists of the modular curves $X_0(N)$ which violate HP, effectively or 
otherwise.  This is because we used Faltings' theorem to deduce the finiteness of $X_0(N)(k)$ for any $k$, 
so long as $g(X_0(N)) \geq 2$.  As mentioned above, the closest function field analogue of Faltings' theorem is 
a finiteness theorem of P. Samuel.  But Samuel's theorem does not apply here because the curves $X_0(N)_{/k}$ 
have algebraic moduli (indeed, they have canonical $\F_p$-models).  
\\ \indent
On the other hand, we can still get Hasse Principle violations using the modular curves $X(N,M)$ as above 
provided that we have a suitable analogue of Merel's theorem, i.e., a result which tells us that for 
fixed $k$ and sufficiently large $M$, the only $k$-rational points on $X_1(M)$ are the cusps.  Fortunately 
for us, such a result has recently been proven by B. Poonen.  We set the stage by introducing some notation:
\\ \\
Let $K$ be a field of characteristic $p > 0$, and let $E_{/K}$ an elliptic curve with transcendental 
$j$-invariant.  Thus $E$ is ordinary.  Let 
\[\rho_E: \Gal_K \ra \Z_p^{\times} \times \prod_{\ell \neq p} GL_2(\Z_{\ell}) = 
\Aut(E(K^{\sep}[\tors]) = \Aut \left(\Q_p/\Z_p \oplus \bigoplus_{\ell \neq p} (\Q_{\ell}/\Z_{\ell})^2 \right)\]
be the homomorphism attached to the $\Gal_K$-module $E(K^{\sep})[\tors]$.  Put 
\[S = \Z_p^{\times} \times \prod_{\ell \neq p} 
SL_2(\Z_{\ell}). \]

\begin{thm}(Poonen, \cite{Poonen07})
\label{PBT}
For every $d \in \Z^+$, there exists a constant $N(p,d)$ such that: for any field 
$k$ of characteristic $p > 0$, any field extension $K/k(t)$ of degree at most $d$, and any elliptic curve 
$E_{/K}$ with transcendental $j$-invariant, the index $[S:\rho_E(\Gal(K^{\sep}/K)) \cap S]$ is at most 
$N(p,d)$.  
\end{thm}
\noindent
From this we readily deduce the desired boundedness result.
\begin{cor}
\label{CHARPBOUNDCOR}
Let $p$ be a prime number and $k_0 = \F_p(t)$.  Then, for every positive integer $d$, there exists a 
constant $B(p,d)$ such that: for any finite field extension $K/k_0$ of degree at most $d$ 
and any prime number $M > \max(3,p,B(p,d))$, the only $K$-rational points on $X_1(M)$ are cusps.
\end{cor}
\noindent
\begin{proof}
Let $K/k_0$ be an extension of degree $d' \leq d$, let $M > \max{3,p,N(p,d)}$ be a prime number, and suppose 
that $X_1(N)$ has a noncuspidal $K$-rational point.  Since $M \geq 4$, 
$Y_1(M)$ is a fine moduli space, and there is a corresponding pair $(E,P)_{/K}$, where $E$ is an elliptic curve 
and $P \in E(K)$ is a point of order $M$.  There is then is a basis 
for $E[M](K^{\sep})$ with respect to which the mod $M$ Galois representation has image contained 
in the subgroup $T(M) \subset \GL_2(\Z/M\Z)$ of matrices 
\[ \{ \left[ \begin{array}{cc} 1 & b \\ 0 & d \end{array} \right] \ | \ b \in \Z/M\Z, \
 d \in (\Z/M\Z)^{\times} \} .\]
Since $\# T(M) \cap \SL_2(\Z/M\Z) = M$, 
the index of the full Galois representation is at least $\# SL_2(\Z/M\Z) / M = M^2-1 > M$.   Applying 
Theorem \ref{PBT}, we conclude that $j(E)$ must be algebraic.  Let $\F$ be the algebraic closure of $\F_p$ 
in $K$.  Because the order of the automorphism group of any elliptic curve divides $24$, it follows that 
there is an extension $\F'/\F$ of degree dividing $24$ and an elliptic curve $E'/\F'$ such that 
$E'_{/F' K} \cong E_{/\F' K}$.  Moreover, since $E'$ is defined over $\F'$, all the torsion is defined 
over an algebraic extension of $\F'$, hence the torsion field $\F'(P)$ is linearly disjoint from 
$\F' K$ over $F$, so $[\F'(P):\F''] = [\F'K(P):\F'K] = 1$, i.e., $P$ is already $\F'$-rational.  
On the other hand, the degree of $\F'$ over $\F_p$ is at most $24[K:k_0] \leq 24d$, so that $\# \F \leq 
p^{24 d}$ and hence (applying the Hasse bound)
\[M \leq \# E(\F') \leq (p^{12d} + 1)^2. \]
Thus it suffices to take 
\[B(p,d) = \max \left(3,p, N(p,d), (p^{12d} + 1)^2 \right). \]
\end{proof}
\noindent
We can now establish the function field case of our Main Theorem by an argument very similar to that used 
in the number field case, replacing our appeal to Merel's Theorem \ref{MERELTHM} 
with an appeal to Corollary \ref{CHARPBOUNDCOR}.  
\\ \\
Indeed, let $k/\F_p(t)$ be a finite separable extension; put $d = [k:\F_p(t)]$.  Let $M$ be a prime 
number which is greater than the constant $B(p,2d)$ of Corollary \ref{CHARPBOUNDCOR}.  By Dirichlet's theorem 
on primes in arithmetic progressions combined with the Chinese Remainder Theorem, there exist infinitely many 
prime numbers $N$ such that 
$N \equiv -1 \pmod M$ and $-N$ is a quadratic residue mod $p$.  If $p = 2$ the last condition is vacuous, 
so we require moreover that $-N \equiv 1 \pmod 8$.  Thus $p$ splits in $\Q(\sqrt{-N})$.  We claim that if 
$N$ is sufficiently large, $X_0(N)$ -- and, \emph{a fortiori}, $X(N,M)$ -- has no $k$-rational $w_N$-fixed 
points.  Let $d(N,p)$ be the least degree of a $w_N$-fixed point on $X_0(N)_{/\F_p}$.  By the modular interpretation 
of Atkin-Lehner fixed points recalled above, we have 
\[d(N,p) = [\F_p(j(E)):\F_p], \] where 
$E/_{\overline{\F_p}}$ is any elliptic curve with complex multiplication by the maximal order of $\Q(\sqrt{-N})$ 
and $j(E)$ is its $j$-invariant.  Now we recall Deuring's correspondence: reduction modulo $p$ induces 
a bijection from the set of the $j$-invariants of $K$-CM elliptic curves in characteristic $0$ for which $p$ 
splits in the CM field $K$ to the set of $j$-invariants of ordinary elliptic curves over $\overline{\F_p}$ 
\cite[Thm. 13.13]{Lang}.  
Since there are, of course, only finitely many $j$-invariants lying in an extension of $\F_p$ of degree at most 
$d$, this shows that for fixed $p$ and $N$ lying in the above infinite set of primes, $d(N,p)$ approaches 
infinity with $N$.  Thus for all sufficiently large $N$, $w_N$ has no $k$-rational fixed points: hypothesis 
$(i)_k$ holds.  Moreover, our choice of $M$ implies that $X(N,M)/w_N$ has only cuspidal $k$-rational points, so 
hypothesis $(iv)_k$ holds.  We have at least one $\F_p$-rational cusp, so hypothesis (iii) holds.  
Finally, hypothesis (ii) certainly holds: we can either use the same argument as we did in 
characteristic $0$ or simply observe that it follows from the characteristic $0$ case by specialization.  
Therefore Corollary \ref{TWISTAHPCOR} applies, completing the proof of the Main Theorem in the function field case.

\subsection{A Complement: Drinfeld modular curves} \textbf{} \\ \\ \noindent
Let $k$ be a global function field of characteristic $p > 2$.  As we saw, Samuel's theorem does not apply 
to the modular curves $X(N,M)_{/k}$.  One may view this as a hint that we have chosen ``the wrong modular curves'' 
for the function-field version of the argument.  Indeed, a family of curves with properties highly analogous to 
those of modular curves over number fields and to which Samuel's theorem does apply are the Drinfeld modular curves.
\\ \\
Specifically, take $A = \F_p[t]$ and for $\mathfrak{n}$ a nonzero prime ideal of $A$ (which we identify 
with the corresponding monic irreducible polynomial) there is a curve 
$X_0(\mathfrak{n})_{/k_0}$ and a canonical $k_0$-rational involutory automorphism $w_{\nn}$.  The involution 
$w_{\nn}$ always has geometric fixed points, and the least degree of a fixed point is equal to the class number 
$h(\nn)$ of the ``imaginary quadratic field'' $k_0(\sqrt{\nn})$.  Both $h(\nn)$ and the genus of 
$X_0(\nn)/w_{\nn}$ tend to infinity with the norm of $\nn$.  Moreover $X_0(\mathfrak{n})$ has at least one $k_0$-rational 
cusp.   (For all these facts, see \cite{Gekeler86}, \cite{Gekeler01}.) Therefore:

\begin{thm}
\label{DRINFELDTHM}
For each fixed odd prime $p$ and finite separable extension $k/k_0$, there exists $N = N(p,k)$ such that 
for all primes $\nn$ of norm larger than $N$, there are infinitely many separable quadratic extensions $l/k_0$ 
such that the Atkin-Lehner twist of $X_0(\nn)$ by $w_{\nn}$ and $kl/k$ violates the Hasse Principle.
\end{thm}
\noindent
However, there are two differences from the classical modular case. \\ \indent
One the one hand, the analogue of Merel's theorem is not yet available: the uniform boundedness of 
torsion points on rank $2$ Drinfeld modules is believed to be true but remains open.
Thus, consideration of Drinfeld modular curves with 
$\Gamma_0(\nn) \cap \Gamma_1(\mathfrak{m})$-level structure would not help to make the constructive effective.
\\ \indent
On the other hand, this help is not needed, because work of Szpiro \cite{Szpiro} makes Samuel's theorem effective (c.f. 
\cite[Lemma 2.1]{Poonen09}).
\\ \\
We therefore have another way of answering Poonen's question for function fields of positive, odd characteristic.

\subsection{Poonen's theorem} \textbf{} \\ \\
The main result of Poonen's paper \cite[Thm. 1.1]{Poonen09} is stronger than the one mentioned in the introduction.  
He proves:
\begin{thm}(Poonen)
\label{POOTHM}
There is an algorithm which takes as input a global field $k$ and $n \in \N$ and returns 
a nice algebraic curve $C_{/k}$ with $C(\mathbb{A}_k) \neq \varnothing$ and $\# C(k) = n$.
\end{thm}
\noindent
The proof of our Main Theorem leads to a proof of Theorem \ref{POOTHM} as well.  This is done in much the same 
way as in \cite{Poonen09}, and we make no claims to novelty.  Nevertheless 
we feel that the result is interesting enough to be included here as well as in \emph{loc. cit.}
\\ \\
The main idea is as follows: we have already proved Theorem \ref{POOTHM} in the case $n = 0$.  In the remaining 
case $n > 0$ -- in which case the hypothesis $C(\mathbb{A}_k) \neq \varnothing$ is of course automatic -- it 
suffices to find the following: for any global field $k$, a curve $C_{/k}$ in which the set $C(k)$ is nonempty, 
finite and effectively computable.  Then by a simple weak approximation argument, we can build a finite, separable 
branched covering $\pi: C' \ra C$ which has prescribed local behavior at each of the finitely many fibers containing a 
$k$-rational point.  In particular, for each $P \in C(k)$, we may arrange for the fiber over $P$ to contain 
$n_P$ rational points, for any $0 \neq n_P \leq \deg(\pi)$.  In particular, for any $n \geq 2$, by building a covering of degree $n$, 
requiring one of the $k$-rational points to split completely into $n$ rational points, and requiring the fibers 
to be irreducible over the remaining $k$-rational points, we get a curve $C'$ with exactly $n$ rational points.  
We proceed similarly for $n = 1$, except we need to use a covering of degree at least $2$.  
\\ \\
Starting from scratch, it is no trivial matter to produce such a curve $C$ over an arbitrary global field.  Fortunately 
we have already done so: we observed above that for fixed $k$ and sufficiently large $N$, the only rational points on 
the curves $X(N,M)_{/k}$ are the cusps, which are finite in number, have effectively computable (indeed, well-known) 
fields of definition, and for which at least one is rational over $k_0$.  This is enough to prove 
Theorem \ref{POOTHM}.
\\ \\
Poonen's construction is in a way not so different from ours: he uses modular curves $X_1(p^n)$.  However, by virtue of 
choosing a large power of a fixed prime he is able to use the results of Dem'janenko-Manin, which are 
considerably more elementary than those of Merel.  It would certainly be of interest to see further examples.

\section{HP-violations with prescribed genus} 
\noindent

\begin{thm}
\label{BIELLIPTICTHM} Let $k$ be a global field of characteristic different from $2$.  
We \textbf{suppose} that there exists an elliptic curve $E_{/k}$ with 
$E(k)$ finite.  There is an algorithm which takes as input an integer $g \geq 2$ and outputs a nice 
curve $C_{/k}$ of genus $g$ which violates the Hasse principle.
\end{thm}
\begin{proof}
Let $E_{/k}$ be as in the statement of the theorem.  Let $q: E \ra \PP^1$ be the degree $2$ map obtained 
by quotienting out by the involution $[-1]$.  (In more elementary terms, $q: (x,y) \mapsto x$.)  By hypothesis, 
for all but finitely many points $x \in \PP^1(K)$, the pullback of the degree one divisor $[x]$ on $\PP^1$ 
is the degree $2$ divisor $[P] + [-P]$ such that the field of definition of $P$ is a quadratic extension of 
$k$.  Moreover, for every finite extension $l/k$, the set $q(E(l))$ of $x$-coordinates of $l$-rational points 
of $E$ is \textbf{thin} in the sense of \cite[$\S 3.1$]{TGT}.  Since $k$ is a Hilbertian field, it follows 
that we can choose $g-1$ points $x_i \in k(\PP^1)$ such that the fields of definition of the points $[P_i]$, $[-P_i]$ 
in the pulled back divisor $D_i = [P_i] + [-P_i] = q^*([x_i])$ are distinct quadratic extensions $l_1,\ldots,l_{g-1}$ 
of $k$.  In particular, the $D_i$'s have pairwise disjoint supports.  \\ \indent
Next, note that all of the $D_i$'s are linearly equivalent to each other, and also to $q^{-1}([\infty]) = 
2[O]$.  Therefore the divisor 
\[D := \sum_{i=1}^{g-1} D_i - (2g-2) [O] \]
is linearly equivalent to $0$ on $E$, so it is the divisor of a function $f \in k(E)$.  Moreover, we have 
the leeway of multiplying $f$ by any nonzero element $a$ of $k$ without changing its divisor, and it is easy to see that by an appropriate 
choice of $a$, we can arrange for the pullback of the divisor $[O]$ in the twofold cover $C$ of $E$ defined 
by the equation $y^2 = \sqrt{a f}$ to consist of two distinct, $k$-rational points.  The curve $C$ is ramified at the $2(g-1)$ distinct points comprising the support of $D_i$, so by Riemann-Hurwitz has genus 
$(g-1) + 1 = g$.  We may view the twofold covering $C \ra E$ as being the quotient by an involution $\iota$ 
on $C$, and one immediately verifies that the pair $(C,\iota)$ satisfies hypotheses (i) through (iv) of Theorem 
\ref{TAHP}.  
\end{proof}

\begin{cor}
\label{CMCOR}
For every number field $k$ and any integer $g \geq 2$, there exists a bielliptic curve $C_{/k}$ violating HP.
\end{cor}
\begin{proof}
Indeed, Mazur and Rubin have recently shown that for every number field $k$, there exist infinitely many 
elliptic curves $E_{/k}$ with $E(k) = 0$ \cite[Cor. 1.9]{MR}.
\end{proof}
\noindent
Remark: Genus one counterexamples to HP over $\Q$ were first constructed by Lind and Reichardt.  
(See $\S 5$ for some further discussion of the genus one case.)  D. Coray and C. Manoil showed \cite[Prop. 4.2, 4.4]{CM} that for every 
$g \geq 2$ there exists a hyperelliptic curve $C_{/\Q}$ of genus $g$ violating HP.  Using Wiles' theorem on 
rational points on Fermat curves, they also gave \cite[Prop. 4.5]{CM} examples of nonhyperelliptic HP violations of genus $4k^2$ 
for any $k \geq 1$.  
\\ \\
Remark: Any nice curve of genus $2$ is hyperelliptic.  A nice nonhyperelliptic curve of genus $3$ is (canonically embedded 
as) a plane quartic curve.  Bremner, Lewis and Morton have shown \cite{BLM} that the nice plane quartic
\[C_{/\Q}: 3X^4 + 4Y^4 = 19Z^4 \]
violates HP.  For $g \geq 4$ a bielliptic curve of genus $g$ is not hyperelliptic \cite{ACGH}, so it follows 
from Corollary \ref{CMCOR} that nonhyperelliptic HP violations over $\Q$ exist in every conceivable genus.
\\ \\
Over any global field $k$ with $\car(k) \neq 2$, we can prove a slightly weaker result.
\begin{thm}
Let $k$ be a global field of odd characteristic and $g \geq 4$ be an integer.  Then there is a curve 
$C_{/k_0}$ of genus $g$ which violates the Hasse Principle over $k$.  
\end{thm}

\begin{proof}
Fix $g_0 \geq 2$, and let $P(x) \in k[x]$ be a monic polynomial of degree $2g_0+1$ with distinct roots in 
$\overline{k}$.  When $\car(k_0) = p > 0$, we require that $P$ is not $PGL_2$-equivalent to a polynomial defined 
over $\overline{\F_p}$.  Then 
\[X: y^2 = P(x) \]
defines a nonisotrivial hyperelliptic curve of genus $g_0$ with a $k$-rational Weierstrass point $O$.  Let $q: X \ra \PP^1$ 
denote the hyperelliptic involution.  Exactly as in the previous construction, for any $d \geq 1$ we may choose 
$d$ points $x_i \in \PP^1(k)$ such that $D_i := q^*([x_i]) = [P_i] + [q(P_i)]$, such that the residue fields 
$l_i$ of $P_i$ are distinct quadratic extensions of $k$.  Now take 
\[D = \sum_{i=1}^d D_i - 2d [O], \]
and choose $f \in k(X)$ with $\div(f) = D$ such that $O$ splits into $2$ $k$-rational points on the corresponding 
$2$-fold cover $C \ra X$ defined by the extension $k(X)(\sqrt{f}) / k(X)$.  By Riemann-Hurwitz, the genus of $C$ 
is $g = 2g_0 + d - 1$.  Writing $\iota$ for 
the induced involution on $C$, the pair $(C,\iota)$ satisfies all the hypotheses of Theorem \ref{TAHP}, so that 
there are infinitely many twists which violate HP over $k$.  Taking $g_0 = 2$, we get curves of every 
genus $g \geq 4$.
\\ \indent
Now if we require that $P$ has $k_0$-coefficients, then the construction can be carried out over $k_0$ by taking $k_0$-rational points on 
$\PP^1$ and ensuring that at least one of the quadratic extensions $l_i/k_0$ is not contained in $k$ (this is 
why we showed in the proof of Theorem \ref{BIELLIPTICTHM} that infinitely many quadratic extensions were 
possible).  Then the hypotheses of Corollary \ref{TWISTAHPCOR} are satisfied, so that we get infinitely many curves of any given $g \geq 4$ which 
are defined over $k_0$ and violate HP over $k$.
\end{proof}

\section{Conjectures on HP violations in genus one}
\noindent
We wish to consider some conjectures involving HP violations on curves of genus one and the logical relations 
between them.
\begin{conj} 
\label{CONJ1} For every global field $k$, there exists a genus one curve $C_{/k}$ which violates the Hasse Principle.  
Equivalently, there exists an elliptic curve $E_{/k}$ with $\Sha(k,E) \neq 0$.
\end{conj}
\noindent
\noindent
The most important conjecture on Shafarevich-Tate groups is, of course, that they are all finite.  More precisely, 
say that a global field $k$ is $\Sha$-\textbf{finite} if $\Sha(k,E)$ is finite for all elliptic curves $E_{/k}$.  It 
is expected that every global field is $\Sha$-finite.  The following observation is mostly for amusement:
\begin{prop}
\label{AMUSEPROP}
Suppose that the prime global subfield $k_0$ is \emph{not} $\Sha$-finite.  Then for every finite field extension 
$k/k_0$, there exists a genus one curve $C_{/k}$ violating the Hasse Principle.
\end{prop}
\begin{proof}
Our assumption is that there exists an elliptic curve $E_{/k_0}$ with $\# \Sha(k_0,E) = \infty$.  Since for all 
$n \in \Z^+$, $\Sha(k_0,E)[n]$ is finite, the infinitude of $\Sha(k_0,E)$ means that it contains elements of arbitarily 
large order.  In particular, if $k/k_0$ is a field extension of degree $d$, then let $C$ be a genus one curve 
corresponding to a class $\eta \in \Sha(k_0,E)$ of order strictly greater than $d$.  Then the index of $C$ is strictly 
greater than $d$, and since $C$ has index one, this implies that the least degree of a closed point exceeds $d$, so 
that we must have $C(k) = \varnothing$.
\end{proof}
\noindent
The proof of Proposition \ref{AMUSEPROP} also shows that Conjecture \ref{CONJ1} is implied by the following
\begin{conj}
\label{CONJ2}
For $n,k \in \Z^+$ with $n > 1$, there exists an elliptic curve $E_{/k_0}$ such that $\Sha(k_0,E)$ has at least 
$k$ elements of order $n$.
\end{conj}
\noindent
Conjecture \ref{CONJ2} however seems to lie much deeper than Conjecture \ref{CONJ1}.  When $k_0 = \Q$, 
Conjecture \ref{CONJ2} is known to hold for $n = 2,3,5$; 
Donnelly and Matsuno \cite{Matsuno} have both announced proofs for $n = 7$; and Matsuno has announced a proof 
for $n = 13$ (\emph{loc. cit.}).  All of these 
arguments use the fact that the modular curve $X_0(n)$ has genus zero for these values of $n$, and, 
unfortunately, these are the only prime numbers $n$ such that $X_0(n)$ has genus $0$.  Matsuno also refers to 
as yet unpublished work of Naganuma claiming similar results 
over many quadratic extensions of $k_0$ when $X_0(n)$ has genus one.  Decidedly new ideas seem to be needed to handle 
the remaining cases.  
\\ \\
Remark: If there is any work on the $k_0 = \F_p(t)$ analogues of these results, I am not aware of it.  Nevertheless 
one expects that similar results can be proven, at least as long as one avoids the case $p \ | \ n$.  
\\ \\
The HP-violations constructed in $\S 3$ all had the property that they were base extensions of HP violations of the 
prime global subfield.  It is natural to consider such constructions in genus one:
\begin{conj}
\label{CONJ3}
For every global field $k$, there is a genus one curve $C$ defined over $k_0$ such that 
$C$ violates the Hasse Principle over $k_0$ and also over $k$.
\end{conj}
\noindent
Evidently Conjecture \ref{CONJ2} $\implies$ Conjecture \ref{CONJ3} $\implies$ Conjecture \ref{CONJ1}.  One pleasant aspect 
of Conjecture \ref{CONJ3} is that it suffices to prove it for $k$ sufficiently large, e.g. Galois over $k_0$ and containing 
sufficiently many roots of unity for Kummer theory to apply. 
\section{Passage to a subvariety}
\noindent

\subsection{Restriction of scalars and HP violations} \textbf{} \\ \\ 
\noindent
Let $k$ be a global field, $l/k$ a finite separable field extension, and $V_{/l}$ a nice variety.  
Recall that there is a $k$-variety $W = (\Res_{l/k} V)$, called the \textbf{Weil restriction} of $V$, whose 
functor of points on affine $k$-schemes is as follows: $\Spec A \mapsto 
V(A \otimes_k l)$.  Taking $A = k$, we get 
\[W(k) = V(k \otimes_k l) = V(l), \]
and taking $A = \mathbb{A}_k$, we get 
\[W(\mathbb{A}_k) = V(\mathbb{A}_k \otimes l) = V(\mathbb{A}_l), \]
so that $V$ violates the Hasse Principle over $l$ iff $W = \Res_{l/k} V$ violates the Hasse Principle 
over $k$.  Moreover, if $V$ is smooth, projective and geometrically integral, then 
then $W$ is smooth, projective and geometrically integral, and of dimensoin $[l:k] \dim V$.  That is to 
say, a Hasse Principle violation over any finite separable extension $l/k$ gives rise a Hasse Principle 
violation over $k$, at the expense of multiplying the dimension by $[l:k]$.

\subsection{HP violations from torsors under abelian varieties}
\begin{thm}
Let $k$ be any global field.  Then there exists an abelian variety 
$A_{/k}$ and a torsor $W_{/k}$ under $A$ such that $W$ violates the Hasse Principle over $k$.
\end{thm}
\begin{proof} 
Suppose $k$ has characteristic exponent $p$, and put $d = [k:k_0]$.  Let $d' > 1$ be any integer which is 
prime to $pd$, and let $E_{/k_0}$ be any elliptic curve.  By \cite[Thm. 3]{CS}, there exists a field 
extension $l_0/k_0$ of degree $d'$ and a genus one curve $C_{/l_0}$ with Jacobian elliptic curve 
$E_{/l_0}$, such that $C$ is an element of $\Sha(l_0,E)$ and has exact order $d'$.  Let 
$A = \Res_{l_0/k_0} E$ and $W = \Res_{l_0/k_0} C$.  Then $A_{/k_0}$ is an abelian variety of dimension 
$d_0$, and $W$ is a torsor under $A$ whose corresponding element of $H^1(k_0,A)$ is locally trivial 
and of order $d'$.  Finally, consider $W_{/l}$: certainly it is locally trivial.  On the other hand 
since $[l:k]$ is prime to $d'$, the cohomological restriction map $H^1(A,k_0)[d'] \ra H^1(A,l)[d']$ is 
injective, so that $W(l) = \varnothing$. 
\end{proof}
\noindent

\subsection{Subvarieties of Varieties Violating HP} \textbf{} \\ \\
\noindent
The construction of $\S 6.2$ can be made effective.  The drawback, 
however, is that $\dim W$ depends on $[l:k_0]$ and approaches infinity as $[l:k_0]$ becomes divisible 
by all sufficiently small primes.  
\\ \\
This raises the following question:
\begin{ques} Let $k$ be a global field and $W_{/k}$ a nice variety which violates the Hasse Principle.  
Must there be a nice curve $C \subset_k W$ which violates the Hasse Principle?
\end{ques}
\noindent
Evidently any subvariety $W \subset_k V$ has $W(k) = \emptyset$, so the issue is one of constructing a 
curve which has points everywhere locally.  
\\ \\
Here is a partial result in this direction:
\begin{thm}
\label{SURFACETHM}
Let $k$ be a number field and $W_{/k}$ a nice variety which violates HP, with $\dim W \geq 2$.  Then 
there is a nice surface $S \subset_k W$ which violates HP.
\end{thm}
\noindent
We require the following lemma, whose proof was supplied by E. Izadi.%
\begin{lemma}
\label{IZADILEMMA}
Let $k$ be a field of characteristic $0$, and let $W_{/k}$ a nice variety of dimension $d$. 
Let $C \subset W$ be an integral curve which is a $k$-subscheme of $W$, and let $Z = \sum_i [P_i]$ 
be an effective $k$-rational zero-cycle with all the $P_i$'s distinct.  Then there exists a closed 
subscheme $S \subset W$ which is a nice surface, and which contains both $C$ and the support of $Z$.
\end{lemma}
\begin{proof} Choose an ample linear system, say $L$.  Then there is a high multiple of $L$, say $M$, 
such that $C$ and $Z$ are cut out (as schemes) by the members of $M$ that contain them. By taking $M$ 
sufficiently large (as a multiple), if we choose $d-2$ general members of $M$ that contain $C$ and $Z$, then 
their intersection is a smooth surface containing $C$ and $Z$.
\end{proof}
\noindent
\emph{Proof of Theorem \ref{SURFACETHM}.} By Bertini, there is a nice curve $C \subset_k W$.  Like any 
geometrically irreducible variety, $C$ can fail to have local points only at a finite set, say $S$, of places of $k$.  
Since $W$ has points everywhere locally, there exists an effective zero-cycle $Z_{/k}$ such that 
for all $v \in S$, $Z(k_v) \neq \varnothing$.  Now apply Lemma \ref{IZADILEMMA}. \qed

\end{document}